# Structural Minimum Controllability Problem for Switched Linear Continuous-Time Systems

Sérgio Pequito, *Member, IEEE,* and George J. Pappas, *Fellow, IEEE*

## Abstract

This paper addresses a structural design problem in control systems, and explicitly takes into consideration the possible application to large-scale systems. More precisely, we aim to determine and characterize the minimum number of manipulated state variables ensuring structural controllability of switched linear continuous-time systems. Towards this goal, we provide a new necessary and sufficient condition that leverages both graph-theoretic and algebraic properties required to ensure feasibility of the solutions. With this new condition, we show that a solution can be determined by an efficient procedure, i.e., polynomial in the number of state variables. In addition, we also discuss the switching signal properties that ensure structural controllability and the computational complexity of determining these sequences. In particular, we show that determining the minimum number of modes that a switching signal requires to ensure structural controllability is NP-hard.

## I. Introduction

Switched systems have been intensively studied, and the primary motivation comes partly from the fact that these systems have numerous applications in the control of mechanical systems, process control, automotive industry, power systems, aircraft and traffic control, and many other fields [13], [27]. Among switched systems, those with all subsystems described by linear differential equations, referred to as *switched linear systems*, have attracted most of the attention [13]. Recent efforts aimed to analyze controllability and reachability properties of these systems [2], [11], [27], [28]. Nonetheless, only recently controllability was studied for the class of uncertain switched linear system, i.e., the parameters of subsystems' state

Both authors are with the Department of Electrical and Systems Engineering, School of Engineering and Applied Science, University of Pennsylvania, (e-mails: {sergo,pappasg}@seas.upenn.edu).

This work was supported in part by the TerraSwarm Research Center, one of six centers supported by the STARnet phase of the Focus Center Research Program (FCRP) a Semiconductor Research Corporation program sponsored by MARCO and DARPA, and the NSF ECCS-1306128 grant.

                                                                 

matrices are either unknown or zero [14]. This assumption copes with scenarios where the system parameters are difficult to identify and obtained with a certain approximation error. Thus, structural properties that are independent of a specific value of unknown parameters are of particular interest. Subsequently, a switched linear system is said to be structurally controllable if one can find a set of values for the unknown parameters such that the corresponding switched linear system is controllable in the classical sense [14]. Whereas in [14] necessary and sufficient conditions to characterize the structural controllability were provided, due to economic constraints we aim to determine and characterize the smallest subset of actuated state variables yielding structural controllability.

The problem of determining the smallest subset of state variables yielding controllability is commonly referred to as *minimum controllability problem*, and it was studied in the context of linear time-invariant systems in [15], [25]. Recently, the minimum controllability problem has also been explored under additional energy constraints, in particular, metrics that depend on the controllability Grammian, see [3], [4], [12], [17], [26], [29]. Alternatively, if we aim to achieve structural controllability, we refer to the problem as *structural* minimum controllability problem. This problem has been fully addressed in the context of linear-time invariant systems in [19] for homogeneous costs; and the computational complexity analyzed for several classes of systems in [1]. In the current manuscript, we extend these results to the case of switched linear continuous-time systems. Notice that whereas in [1], [19] graph-theoretic properties of structural controllability rely in directed graphs interpretations of the system, these no longer hold to characterize structural controllability of structural switched linear continuous-time systems, see [14] for details. In particular, the analysis of structural controllability of structural switched linear continuous-time systems cannot be reduced to the analysis of a structured linear system. Consequently, in this paper we provide a systematic approach that leverages the combination of graph-theoretic and algebraic conditions to obtain and characterize the solutions to the structural minimal controllability problem for structural switched linear continuous-time system.

The current work also differs from [23], [24], where the structural minimum controllability problem aimed to ensure structural controllability for each mode of the switched continuous-time linear system; note that this conservative notion contrasts with the controllability definition considered in the present manuscript. In [20] the structural minimal controllability problem for linear time-invariant systems was considered under heterogenous cost, i.e., the variables actuated can incur in different costs. In particular, in [20] this problem is shown to be poly-nomially solvable, and in [16], the computational complexity was improved when binary costs are considered. More recently, in [22] the problem was extended to the case where a state variable has a cost that depends on the input that actuates it, hence, leading to a multiple heterogenous cost scenario. Alternatively, the problem of determining the minimum number of actuators from a given collection of possible actuator-state configurations was shown to be (in general) NP-hard [18]. Notwithstanding, in [21] it was shown that the same problem can be polynomially solvable when the dynamic matrix is irreducible.

The main contributions of this paper are fourfold: (*i*) we provide a new necessary and sufficient condition that leverages both graph-theoretic and algebraic properties required to ensure structural controllability of switching linear continuous-time systems; (*ii*) we characterize the solutions to the structural minimum controllability problem for switched linear continuous-time systems. In particular, we characterize *dedicated solutions*, i.e., an actuator can only actuate a single state variable, and *minimal solutions*, i.e., the minimum number of actuators actuating the minimum number of state variables; (*iii*) we propose an algorithm that leverages both graph-theoretic and algebraic properties of structural controllability of switching linear continuous-time systems to determine a solution in $(mn)^{\alpha}$, where $n$ denotes the number of state variables, $m$ denotes the number of modes of the switching linear continuous-time system, and $\alpha < 2.373$ is the exponent of the $n \times n$ matrix multiplication; and (*iv*) we provide new insights on how the switching sequences affect the controllability of a structural switching linear continuous-time system. In particular, we show that determining the minimum collection of nodes in a sequence of modes ensuring structural controllability is NP-hard.

The rest of the paper is organized as follows: Section II provides the formal statement of the problem addressed in this paper. Next, Section III reviews some concepts, introduces key results in structural systems theory and establishes their relations to graph-theoretic constructs. In Section IV, we present the main results. Next, we present an illustrative example in Section V. Finally, Section VI concludes the paper and discusses avenues for further research.

## II. PROBLEM STATEMENT

In this section, we formally introduce the structural minimum controllability problem for switched linear continuous-time systems.

Consider the following switched linear continuous-time system

$$\dot{x}(t) = A_{\sigma(t)}x(t) + B_{\sigma(t)}u(t), \qquad (1)$$

where $\sigma : \mathbb{R}^{+} \to \mathbb{M} \equiv \{1, \ldots, m\}$ is a switching signal, $x(t) \in \mathbb{R}^{n}$ the state of the system at the instance of time $t$, and $u(t) \in \mathbb{R}^{n}$ represents the piecewise continuous input signal. In the sequel, we identify (1) by the pair $(A_{\sigma(t)}, B_{\sigma(t)})$, that contains $m$ modes with subsystems $(A_i, B_i)$, $i \in \{1, \ldots, m\}$, and $\sigma(t) = i$ implies that the $i$th subsystem $(A_i, B_i)$ is active at time instance $t$. Further, the switched linear continuous-time system (1) is said to be controllable (or equivalently, $(A_{\sigma(t)}, B_{\sigma(t)})$ is controllable) if for any initial state $x(0) = x_0$, and a desired state $x_d$, there exists a time instance $t_f > 0$, a switching signal $\sigma : [0, t_f) \to M$ and an input $u : [0, t_f) \to \mathbb{R}^{p}$ such that $x(t_f) = x_d$. This notion of controllability enables the analysis of switching systems where we either have access to 'common' transitions and knowledge of the existing modes of the switching system, or the cases where the controller is equipped with supervisory capabilities enabling the system to switch between modes.

Due to economic constraints, one is interested in deploying the minimum actuation capabilities that enable the controllability of the system, which can be captured by the following optimization problem. Given the switched linear continuous-time system (1), we aim to determine the sparsest input matrices $\{\bar{B}_i^*\}_{i=1}^{m}$ for the different $m$ modes required to ensure controllability, as a solution to the following optimization problem:

$$\min_{B_1, \ldots, B_m \in \mathbb{R}^{n \times n}} \quad \sum_{i=1}^{m} \|B_i\|_0$$
$$\text{s.t.} \quad (A_{\sigma(k)}, B_{\sigma(k)}) \text{ is controllable,}$$

where $\|M\|_0$ is the zero (quasi) norm, i.e., it counts the number of non-zeros entries in matrix $M$. Unfortunately, this problem is NP-hard even when $m = 1$, see [15] for details.

Furthermore, the parameters associated with the linear time-invariant modes are often accurately known, which motivates the use of structural system theory [8]. Structural linear systems are linear parameterized systems with a given structure, i.e., the entries of the state space matrix are either free parameters or fixed zeros. Let $\bar{A}_{\sigma(t)} \in \{0, 1\}^{n \times n}$ denote the zero/nonzero structure or *structural pattern* of the system matrix $A_{\sigma(t)}$, whereas $\bar{B}_{\sigma(t)} \in \{0, 1\}^{n \times p}$ is the structural pattern of the input matrix $B_{\sigma(t)}$. More precisely, an entry in these matrices is zero if the corresponding entry in the system matrices is equal to zero, and described by an arbitrary parameter (denoted by one) otherwise. Therefore, a pair $(\bar{A}_{\sigma(t)}, \bar{B}_{\sigma(t)})$ is said to be structurally controllable if there exists a pair $(A'_{\sigma(t)}, B'_{\sigma(t)})$ respecting the structure of $(\bar{A}_{\sigma(t)}, \bar{B}_{\sigma(t)})$, i.e., same locations of zeros and nonzeros, such that $(A'_{\sigma(t)}, B'_{\sigma(t)})$ is controllable. Further, it can be shown that if a pair $(A_{\sigma(t)}, B_{\sigma(t)})$ is structurally controllable, then almost all (with respect to the Lebesgue measure) pairs with the same structure as $(\bar{A}_{\sigma(t)}, \bar{B}_{\sigma(t)})$ are controllable [14]. In essence, structural controllability is a property of the structure of the pair $(\bar{A}_{\sigma(t)}, \bar{B}_{\sigma(t)})$ and not of the specific numerical values.

Subsequently, the *structural minimum controllability problem for switched linear continuous-time systems* problem can be stated as follows:

$\mathcal{P}_1$ Given the structure of the matrices of the switched linear system in (1), i.e., $\{\bar{A}_i\}_{i=1}^m$, we aim to determine the sparsest collection of input matrices $\{\bar{B}_i^*\}_{i=1}^m$ required to ensure its structural controllability, i.e., that is the solution to the following problem:

$$\min_{\bar{B}_1,\dots,\bar{B}_m \in \{0,1\}^{n \times n}} \sum_{i=1}^m \|\bar{B}_i\|_0$$
$$\text{s.t.} \quad (\bar{A}_{\sigma(k)}, \bar{B}_{\sigma(k)}) \text{ is struct. controllable.}$$

Notice that a solution to $\mathcal{P}_1$ may consist of $\bar{B}_i$ with columns with all entries are equal to zero, which can be disregarded when considering the deployment of the inputs required to actuate the system. In addition, in the worst case scenario, we obtain structural controllability by taking the identity matrix as the input matrix, which justifies the dimensions chosen for the solution search space. Furthermore, some solutions may comprise at most one nonzero entry in each column; in other words, solutions in which each input actuates at most one state variable. These inputs are referred to as *dedicated inputs*, and they correspond to the columns of the input matrix $\bar{B}_i$ with exactly one nonzero entry. Additionally, if a solution $\{\bar{B}_i^*\}_{i=1}^m$ is such that all its nonzero columns consist of exactly one nonzero entry, it is referred to as a *dedicated solution*, otherwise it is referred to as a *non-dedicated solution*.

Finally, note that the solution procedure for $\mathcal{P}_1$ also addresses the corresponding structural observability output matrix design problem by invoking the duality between observability and controllability in linear time-invariant systems for each mode of the switching linear continuous-time system [10].

## III. PRELIMINARIES AND TERMINOLOGY

In this section, we review some notions of controllability of switched linear continuous-time systems, and their counterpart using structural systems theory [8].

To assess the controllability for switched linear continuous-time systems consider the following definitions.

**Definition 1** ([14])**.** *The* controllability matrix *for switched linear continuous-time system as described in* (1) *is given by*

$$\mathcal{C}(A_{\sigma(k)}, B_{\sigma(k)}) = [B_1, B_2, \ldots, B_m, A_1B_1, A_2B_1, \ldots,$$

$$A_mB_1, A_1B_2, A_2B_2, \ldots, A_mB_2, \ldots, A_1B_m,$$

$$A_2B_m, \ldots, A_mB_m, A_1^2B_1, A_2A_1B_1, \ldots, A_mA_1B_1,$$

$$A_1A_2B_1, A_2^2B_1, \ldots, A_mA_2B_1, \ldots, A_1A_mB_m,$$

$$A_2A_mB_m, \ldots, A_m^2B_m, \ldots, A_1^{n-1}B_1, A_2A_1^{n-2}B_1,$$

$$\ldots, A_1A_2A_1A_1^{n-3}B_1, A_2^2A_1^{n-3}B_1, \ldots, A_1^{n-3}B_1,$$

$$\ldots, A_1A_m^{n-2}B_m, \ldots, A_m^{n-1}B_m].$$

$\diamond$

Additionally, we have the following result.

**Theorem 1** ([14])**.** *The system described by* (1) *is controllable if and only if* rank $\mathcal{C}(A_{\sigma(k)}, B_{\sigma(k)}) = n$. $\diamond$

Now, we associate with the pair $(\bar{A}_i, \bar{B}_i)$, with $\bar{A}_i, \bar{B}_i \in \{0,1\}^{n \times n}$, a directed graph (digraph) $\mathcal{D}(\bar{A}_i, \bar{B}_i) = (\mathcal{V}_i, \mathcal{E}_i)$, referred to as the *system digraph*, with vertex set $\mathcal{V}_i$ and edge set $\mathcal{E}_i$, where $\mathcal{V}_i = \mathcal{U}_i \cup \mathcal{X}_i$ with $\mathcal{X}_i = \{x_1^i, \ldots, x_n^i\}$ and $\mathcal{U}_i = \{u_1^i, \ldots, u_n^i\}$ represents the *state* and *input vertices*, respectively. In addition, $\mathcal{E}_i = \mathcal{E}_{\mathcal{X}_i, \mathcal{X}_i} \cup \mathcal{E}_{\mathcal{U}_i, \mathcal{X}_i}$ where $\mathcal{E}_{\mathcal{X}_i, \mathcal{X}_i} = \{(x_k^i, x_j^i) : [\bar{A}_i]_{jk} \neq 0\}$ and $\mathcal{E}_{\mathcal{U}_i, \mathcal{X}_i} = \{(u_k^i, x_j^i) : [\bar{B}_i]_{jk} \neq 0\}$ represents the *state edges* and *input edges*, respectively. Similarly, we can define a *state digraph* $\mathcal{D}(\bar{A}_i) = (\mathcal{X}_i, \mathcal{E}_{\mathcal{X}_i, \mathcal{X}_i})$. A *directed path* is a sequence of directed edges where every edge ends in a vertex that is the starting of another edge and no vertex is used twice. A state vertex is said to be *non-accessible* by an input vertex if there exists no directed path from an input to the state vertex. Later, given matrices $\bar{M}_1, \ldots, \bar{M}_m \in \{0,1\}^{n \times n}$, their structure can be combined into a single matrix $\bar{M}$ in terms of $\bar{M} = \bar{M}_1 \vee \ldots \vee \bar{M}_m$, where $\vee$ corresponds to the entry-wise operation where if at least one of the entries is non-zero, then it provides a non-zero entry, and zero otherwise. In addition, $[\bar{M}_1, \ldots, \bar{M}_m]$ denotes the concatenation of matrices $\bar{M}_1, \ldots, \bar{M}_m$.

Next, we introduce the notion of a bipartite graph associated with a $m_1 \times m_2$ matrix $\bar{M}$ given by $\mathcal{B}(\bar{M}) = (\mathcal{C}, \mathcal{R}, \mathcal{E}_{\mathcal{C},\mathcal{R}})$, where $\mathcal{R} = \{r_1, \ldots, r_{m_1}\}$ and $\mathcal{C} = \{c_1, \ldots, c_{m_2}\}$ correspond to the labeling row vertices and column vertices, respectively, and $\mathcal{E}_{\mathcal{C},\mathcal{R}} = \{(c_j, r_i) : \bar{M}_{ij} \neq 0\}$. The bipartite graph is an undirected graph with vertex set given by the union of the partition sets $\mathcal{C}$ and $\mathcal{R}$, which we refer to as left and right vertex sets, respectively. A matching $\mathcal{M} \subset \mathcal{E}_{\mathcal{C},\mathcal{R}}$ is collection of edges that have no vertex in common. A maximum matching is a matching with maximum cardinality among all possible matchings. For ease of reference, if a vertex in the left and right vertex set does not belong to an edge in a maximum matching, we then refer to it a right- and left-unmatched vertex, respectively. Additionally, we can consider weights associated with the edges in a bipartite graph to obtain a *weighted bipartite graph* $(\mathcal{B}(\bar{M}) = (\mathcal{C}, \mathcal{R}, \mathcal{E}_{\mathcal{C},\mathcal{R}}), w)$, where $w : E_{\mathcal{C},\mathcal{R}} \to \mathbb{R}$. Subsequently, we can consider the problem of determining the maximum matching with the minimum sum of the weights, that we refer to as the *minimum weight maximum matching*. The minimum weight maximum matching can be generally solvable in $\mathcal{O}(\max\{m_1, m_2\}^\alpha)$, where $\alpha < 2.373$ is the exponent of the $n \times n$ matrix multiplication [9].

In addition, a digraph $\mathcal{D}_S = (\mathcal{V}_S, \mathcal{E}_S)$ is a *subgraph* of $\mathcal{D} = (\mathcal{V}, \mathcal{E})$ if $\mathcal{V}_S \subseteq \mathcal{V}$ and $\mathcal{E}_S \subseteq \mathcal{E}$. Finally, a *strongly connected component* (SCC) is a maximal subgraph (there is no other subgraph, containing it, with the same property) $\mathcal{D}_S = (\mathcal{V}_S, \mathcal{E}_S)$ of $\mathcal{D}$ such that for every $u, v \in \mathcal{V}_S$ there exists a path from $u$ to $v$ and from $v$ to $u$. We can create a *directed acyclic graph* (DAG) by visualizing each SCC as a virtual node, where there is a directed edge between vertices belonging to two SCCs if and only if there exists a directed edge connecting the corresponding SCCs in the digraph $\mathcal{D} = (\mathcal{V}, \mathcal{E})$, the original digraph. The DAG associated with $\mathcal{D}(\bar{A})$ can be computed efficiently in $\mathcal{O}(|\mathcal{V}| + |\mathcal{E}|)$ [6]. The SCCs in the DAG may be further categorized as follows: an SCC is *non-top linked* if it has no incoming edge to its vertices from the vertices of another SCC.

Finally, consider a $m_1 \times m_2$ matrix $\bar{M}$, and let $[\bar{M}] = \{P \in \mathbb{R}^{m_1 \times m_2} : P_{ij} = 0 \text{ if } \bar{M}_{ij} = 0\}$, then the *generic rank* (g-rank) of $\bar{M}$ is given by g-rank$(\bar{M}) = \max\limits_{P \in [\bar{M}]} \text{rank}(P)$.

Now, we revisit necessary and sufficient conditions for the structural controllability of switched linear continuous-time systems.

**Theorem 2** ([14]). *A switched linear continuous-time system* (1) *is structurally controllable if and only if the following two conditions hold:*

(i) $\mathcal{D}(\bar{A}_1 \vee \ldots \vee \bar{A}_m, \bar{B}_1 \vee \ldots \vee \bar{B}_m)$ *has no non-accessible state vertex;*

(ii) *g-rank* $\left([\bar{A}_1, \ldots, \bar{A}_m, \bar{B}_1, \ldots, \bar{B}_m]\right) = n.$   ◇

We note that whereas verifying if the conditions in Theorem 2 hold can be done efficiently [14], designing the sequence of sparsest input matrices such that those conditions yield is a more challenging problem. In fact, a greedy strategy may consist in sequently try to ensure each condition. Nonetheless, optimality of such strategies cannot (in general) be ensured. Therefore, one should resort to such strategies only when the problem at hand is computationally intractable, for instance, NP-hard. In this paper, we will show that $\mathcal{P}_1$ can be polynomially solvable by leveraging both graph-theoretic and algebraic characterizations of the conditions in Theorem 2 captured by the following results.

**Lemma 1** ([7]). *Let* $\bar{M} \in \{0, 1\}^{m_1 \times m_2}$. *There exists a maximum matching of* $\mathcal{B}(\bar{M})$ *with size $n$ if and only if g-rank$(\bar{M}) = n$.*   ◇

**Lemma 2** ([19]). *The digraph* $\mathcal{D}(\bar{A}_1 \vee \ldots \vee \bar{A}_m, \bar{B}_1 \vee \ldots \vee \bar{B}_m)$ *has no non-accessible state vertex if and only if there exits an edge to a state vertex in each non-top linked SCC of the DAG representation of* $\mathcal{D}(\bar{A}_1 \vee \ldots \vee \bar{A}_m)$ *from an input in* $\mathcal{D}(\bar{A}_1 \vee \ldots \vee \bar{A}_m, \bar{B}_1 \vee \ldots \vee \bar{B}_m)$.   ◇

## IV. MAIN RESULTS

In this section, we present the main results of this paper. More specifically, we characterize all the solutions to $\mathcal{P}_1$. This goal is achieved in three steps. First, we determine a dedicated solution that enables structural controllability by performing actuation into a single mode, i.e., we determine a dedicated solution $\bar{B}^*$ such that $\bar{B}_1^* = \bar{B}$ and $\bar{B}_i^* = 0$ for $i = 2, \ldots, m$ ensures structural controllability of $(\bar{A}_{\sigma(t)}, \bar{B}_{\sigma(t)})$ (see Algorithm 1 whose correctness is provided in Theorem 4). Second, we describe in Theorem 5 the non-dedicated solutions $\bar{B}_1^*$, and $\bar{B}_i^* = 0$ for $i = 2, \ldots, m$, which can be obtained from the dedicated solutions $\bar{B}^*$. In the last step, in Theorem 6, we consider the former characterization to describe all possible solutions to $\mathcal{P}_1$. Finally, we present a discussion of results regarding the switching signals required to ensure structural controllability of $(\bar{A}_{\sigma(t)}, \bar{B}_{\sigma(t)}^*)$. In particular, we show that determining the minimum sequence of modes in such switching is NP-hard.

We start by leveraging graph-theoretic and algebraic conditions presented in Lemma 1 and Lemma 2, to rewrite Theorem 2 as follows.

**Theorem 3.** *A switched linear continuous-time system* (1) *is structurally controllable if and only if the following two conditions hold:*

(i) *there exits an edge to a state vertex in each non-top linked SCC of the DAG representation of $\mathcal{D}(\bar{A}_1 \vee \ldots \vee \bar{A}_m)$ from an input in $\mathcal{D}(\bar{A}_1 \vee \ldots \vee \bar{A}_m, \bar{B}_1 \vee \ldots \vee \bar{B}_m)$;*

(ii) *there exists a maximum matching of $\mathcal{B}([\bar{A}_1, \ldots, \bar{A}_m, \bar{B}_1, \ldots, \bar{B}_m])$ with size $n$.* ◇

The conditions provided in Theorem 3 can now be used to obtain a dedicated solution to $\mathcal{P}_1$. More specifically, we propose Algorithm 1 to obtain $B$ such that $\bar{B}_1 = \bar{B}$ and $\bar{B}_i = 0$ for all $i = 2, \ldots, m$ is a dedicated solution to $\mathcal{P}_1$. Towards this goal, Algorithm 1 consists of a two-step algorithm that determines the smallest collection of state variables required to ensure both conditions in Theorem 3. In the first step, it find the set containing the maximum collection of state variables that simultaneously contribute to satisfy both conditions in Theorem 3. We show that this set can be obtained by considering a MWMM on a weighted bipartite graph. The weighted bipartite graph requires the careful crafting of a 'surrogate' matrix $\bar{S}$ that will encode the graph-theoretic properties required to ensure condition (i) in Theorem 3, and a weight function enables the connection of condition (i) and the algebraic condition in Theorem 3-(ii). Consequently, a MWMM determines the maximum set of state variables that need to be actuated, while satisfying both conditions in Theorem 3. More specifically, Algorithm 1 finds a MWMM of a bipartite graph $\mathcal{B}([\bar{A}_1, \ldots, \bar{A}_m, \bar{S}])$, where the matrix $\bar{S}$ has as many columns as the number of non-top linked SCCs in the DAG representation of $\mathcal{D}(\bar{A}_1 \vee \ldots \vee \bar{A}_m)$. The non-zero entries in column $i$ of $\bar{S}$ correspond to the indices of the state variables that belong to the $i$-th non-top linked SCC. In addition, we consider weights in the edges of the bipartite graph: the edges associated with the nonzero entries of $\bar{A}_i$ have zero weight and the edges associated with nonzero entries in $\bar{S}$ have unitary weight. In particular, if an edge in the MWMM contains the vertex corresponding to the columns of $\bar{S}$, then the row vertex $r_i$ in the edge indicates that the state variable $x_i$ contributes to satisfy simultaneously both conditions in Theorem 3. Hence, in the second step, one remains to determine (independently) the smallest sub collection of state variables fulfilling both conditions in Theorem 3.

The next result establishes the correctness and analyzes the implementation complexity of

---

**ALGORITHM 1:**

**Input**: Structural linear time-invariant dynamics in each mode of the structural switching linear continuous-time system described by $\{\bar{A}_i\}_{i=1}^m$.

**Output**: Input matrix describing a dedicated solution $\mathbb{D}(\mathcal{J})$, where $\mathbb{D}(\mathcal{J})$ represents the $n \times n$ diagonal matrix with entries in $\mathcal{J}$ different from zero

**Step 1**. Determine the non-top linked SCCs $\mathcal{N}_i^T$, $i \in \mathcal{I} \equiv \{1, \cdots, \beta\}$, of $\mathcal{D}(\bar{A}_1 \vee \ldots \vee \bar{A}_m) = (\mathcal{X}, \mathcal{E}_{\mathcal{X},\mathcal{X}})$.

**Step 2**. Consider a weighted bipartite graph $\mathcal{B}([\bar{A}_1, \ldots, \bar{A}_m, \bar{S}]) = (\mathcal{C}, \mathcal{R}, \mathcal{E}_{\mathcal{C},\mathcal{R}})$, where $\bar{S}$ is a $n \times \beta$ matrix and $\bar{S}_{i,j} = 1$ if $x_i \in \mathcal{N}_j^T$, and the column vertices be re-labeled as follows: the columns of $\bar{A}_i$ are indexed by $\{c_1^i, \ldots, c_n^i\}$, and the columns of $\bar{S}$ are indexed by $\{s_1, \ldots, s_\beta\}$. In addition, let the weight of the edges $e \in \left( \bigcup\limits_{i=1,\ldots,m} \{c_1^i, \ldots, c_n^i\} \right) \times \mathcal{R}$ be equal to zero, and the weight on the edges $e \in \{s_1, \ldots, s_\beta\} \times \mathcal{R}$ be equal to one.

**Step 3**. Let $\mathcal{M}'$ be the maximum matching incurring in the minimum cost of the weighted bipartite graph presented in Step 2.

**Step 4**. Take $\mathcal{J}' = \{i : (s_j, r_i) \in \mathcal{M}', j \in \{1, \ldots, \beta\}\}$, i.e., the row vertices associated with $S$ that belong to the edges in the MWMM $\mathcal{M}'$ (i.e., those with weight one). In addition, let $\mathcal{J}'' = \{1, \ldots, n\} \setminus \{i \in \{1, \ldots, n\} : (c_j^k, r_i) \in \mathcal{M}', k \in \{1, \ldots, m\}, j \in \{1, \ldots, n\}\}$, and $\mathcal{J}'''$ contains the index of a single state variable from each non-top linked SCC $\mathcal{N}_p^T$, with $p \in \{1, \ldots, \beta\} \setminus \mathcal{J}'$.

**Step 5**. Set $\mathcal{J} = \mathcal{J}' \cup \mathcal{J}'' \cup \mathcal{J}'''$.

---

Algorithm 1.

**Theorem 4.** *Algorithm 1 is correct, i.e., it provides a dedicated solution $\{\bar{B}_i^*\}_{i=1}^m$, with $\bar{B}_1^*$ obtained using Algorithm 1 and $\bar{B}_i^* = 0$ $(i = 2, \ldots, m)$, is a solution to $\mathcal{P}_1$. Furthermore, its computational complexity is $\mathcal{O}((mn+\beta)^\alpha)$, where $\alpha < 2.373$ is the exponent of the $n \times n$ matrix multiplication.* ◇

**Proof:** The correctness of Algorithm 1 follows from the fact that the indices in $\mathcal{J}'$ identify the minimum set of dedicated inputs, simultaneously maximizing the increase in the g-rank of $[\bar{A}_1, \ldots, \bar{A}_m, \mathbb{D}(\mathcal{J}')]$ with respect to $[\bar{A}_1, \ldots, \bar{A}_m]$ by $|\mathcal{J}'|$, and the dedicated inputs assigned to state variables in different non-top linked SCCs. This follows from observing that (by construction) $\mathcal{B}([\bar{A}_1, \ldots, \bar{A}_m])$ results in a minimum weight maximum matching $\mathcal{M}$ with zero weight and size $|\mathcal{M}|$; hence, by Lemma 1, it follows that g-rank$([\bar{A}_1, \ldots, \bar{A}_m]) = |\mathcal{M}|$. Subsequently, a minimum weight maximum matching $\mathcal{M}'$ of $\mathcal{B}([\bar{A}_1, \ldots, \bar{A}_m, \bar{S}])$ equals $|\mathcal{M}'| - |\mathcal{M}|$; hence, increasing by $|\mathcal{M}'| - |\mathcal{M}|$ the g-rank of $[\bar{A}_1, \ldots, \bar{A}_m, \mathbb{D}(\mathcal{J}')]$ with respect to $[\bar{A}_1, \ldots, \bar{A}_m]$, and

contributing to satisfy both conditions in Theorem 3. Nevertheless, it may be insufficient to ensure condition (ii) in Theorem 3, which is fulfilled by taking into account the minimum (additional) collection of dedicated inputs indexed by $\mathcal{J}''$. In addition, by construction of $\bar{S}$ it follows that $\mathbb{D}(\mathcal{J}')$ corresponds to dedicated inputs that are assigned to state variables in different non-top linked SCCs. Hence, $|\mathcal{J}'|$ non-top linked SCCs have incoming edges from different inputs in the system digraph. Thus, contributing to satisfy condition (i) in Theorem 3, but may not be enough to ensure this condition, which is accounted for by considering the minimum collection of dedicated inputs indexed by $\mathcal{J}'''$ that ensures condition (i) in Theorem 3.

In summary, the total number of additional dedicated inputs $\mathbb{D}(\mathcal{J}'')$ required such that g-rank($[\bar{A}_1, \ldots, \bar{A}_m, \mathbb{D}(\mathcal{J}' \cup \mathcal{J}'')]) = n$ is minimized by considering Step 4. Similarly, the total number of additional dedicated inputs $\mathbb{I}_n^{\mathcal{J}'''}$ required such that there exist no non-accessible state vertices in $\mathcal{D}(\bar{A}_1 \vee \ldots \vee \bar{A}_m, \mathbb{D}(\mathcal{J}' \cup \mathcal{J}'''))$ is minimized by considering Step 4. Notice that $\mathbb{D}(\mathcal{J}'')$ are not assigned to non-top linked SCCs previously assigned, otherwise they would have been considered in $\mathbb{D}(\mathcal{J}')$. Finally, by setting $\mathcal{J} = \mathcal{J}' \cup \mathcal{J}'' \cup \mathcal{J}'''$, as in Step 5, $(\mathbb{D}(\mathcal{J}), 0, \ldots, 0)$ is by construction a solution to $\mathcal{P}_1$, since both conditions in Theorem 3 hold. Further, it is minimal since $|\mathcal{J}|$ is minimal, which implies that $\sum_{i=1}^{m} \|\bar{B}_i\|_0 = |\mathcal{J}|$.

The computational complexity follows from noticing that Step 2 can be solved using the Hungarian algorithm that finds a MWMM of $\mathcal{B}([\bar{A}_1, \ldots, \bar{A}_m, \bar{S}])$ in $\mathcal{O}(\max\{|\mathcal{C}|, |\mathcal{R}|\}^\alpha)$, where $|\mathcal{C}|$ denotes the number of column vertices in $\mathcal{B}([\bar{A}_1, \ldots, \bar{A}_m, \bar{S}])$, and $\alpha < 2.373$ is the exponent of the $n \times n$ matrix multiplication, whereas all other steps have linear complexity; hence, Step 2 dominates the final computational complexity, leading to the final complexity of $\mathcal{O}(|\mathcal{C}|^\alpha)$, since $|\mathcal{C}| \geq |\mathcal{R}|$. ∎

**Remark 1.** *We notice that if the structural switching linear continuous-time system only possesses one mode, then it boils down to a structural linear time-invariant, and the characterizations obtained in [19] can be retrieved. Contrarily to the approach presented in [19] that is motivated by the graph-theoretic characterization of the system digraph, Algorithm 1 requires both graph-theoretic and algebraic characterizations, since graph-theoretic properties for structural switching linear continuous-time system are quite diverse from those known for structural linear time-invariant, see [14] for details.* ○

Next, we characterize all the possible sparsest matrices that are solutions to $\mathcal{P}_1$ when a single node is actuated.

**Theorem 5.** *Given $\mathcal{J}', \mathcal{J}''$ and $\mathcal{J}'''$ as in Algorithm 1, then $\bar{B}_1^* = \mathbb{D}(\mathcal{J}' \cup \mathcal{J}'') \vee \mathbb{O}(\mathcal{J}''')$, where $\mathbb{O}(\mathcal{I})$ is the $n \times n$ matrix with exactly one non-zero entry in row indexed in $\mathcal{J}$ and zeros otherwise, and $\bar{B}_i^* = 0$ for $i = 2, \ldots, m$ is a solution to $\mathcal{P}_1$.* ◇

**Proof:** The proof follows by noticing that both conditions of Theorem 3 hold. Condition (ii) in Theorem 3 holds because a maximum matching of $\mathcal{B}([\bar{A}_1, \ldots, \bar{A}_m, \mathbb{D}(\mathcal{J}' \cup \mathcal{J}'')])$, is also a maximum matching of $\mathcal{B}([\bar{A}_1, \ldots, \bar{A}_m, \bar{B}])$. Therefore, since the maximum matching of the former has size $n$, the latter also has size $n$. Secondly, Theorem 3-(i) also holds for $\mathcal{D}(\bar{A}_1 \vee \ldots \vee \bar{A}_m, \bar{B})$; more precisely, $\mathcal{D}(\bar{A}_1 \vee \ldots \vee \bar{A}_m, \mathbb{D}(\mathcal{J}' \cup \mathcal{J}'' \cup \mathcal{J}'''))$ satisfies Theorem 3-(i), see Theorem 4, which implies that there exists a directed path from an input to every state variable. Therefore, by considering $\bar{B} = \mathbb{D}(\mathcal{J}' \cup \mathcal{J}'') \vee \mathbb{O}(\mathcal{J}''')$ there exist an input edge from an input to the same state variables indexed by $\mathcal{J}'''$. Nevertheless, the inputs from where the input edges start are indexed by the columns with non-zero entries in $\mathbb{O}(\mathcal{J}''')$. ∎

Note that in Theorem 5, the matrix $\mathbb{O}(\mathcal{J}''')$ has no constraints on the number of non-zero entries in each column. Hence, $\bar{B} = \mathbb{D}(\mathcal{J}' \cup \mathcal{J}'') \vee \mathbb{O}(\mathcal{J}''')$ is not necessarily a dedicated solution. Subsequently, we can characterize the *minimal* solutions to $\mathcal{P}_1$ as follows.

**Corollary 1.** *Given $\mathcal{J}', \mathcal{J}''$ and $\mathcal{J}'''$ as in Algorithm 1, then $\bar{B}_1^* = \mathbb{D}(\mathcal{J}' \cup \mathcal{J}'') \vee \mathbb{M}(\mathcal{J}''')$, where $\mathbb{M}(\mathcal{I})$ is the $n \times n$ matrix with exactly one non-zero entry in row indexed in $\mathcal{J}$ in some column indexed by $\mathcal{J}' \cup \mathcal{J}''$ and zeros otherwise and $\bar{B}_i^* = 0$ for $i = 2, \ldots, m$ is a minimal solution to $\mathcal{P}_1$.* ◇

**Proof:** The proof follows from Theorem 5, and noticing that the entries associated with columns indexed by $\mathcal{J}' \cup \mathcal{J}''$ cannot be shared by the same column, as consequence of the construction in Theorem 4; in particular, it would compromise condition (ii) in Theorem 3. ∎

Finally, we provide the most general characterization of the sparsest input matrices that are solution to $\mathcal{P}_1$.

**Theorem 6.** *Let $\mathcal{J}', \mathcal{J}''$ and $\mathcal{J}'''$ as in Algorithm 1, and $\bar{B} = \mathbb{D}(\mathcal{J}' \cup \mathcal{J}'') \vee \mathbb{O}(\mathcal{J}''')$, where $\mathbb{O}(\mathcal{I})$ is the $n \times n$ matrix with exactly one non-zero entry in row indexed in $\mathcal{J}$ and zeros otherwise. If $\bar{B}_i^*$, for $i = 1, \ldots, m$, is such that the following holds:*

- *$\bar{B}_i^*$ contains as columns only the columns of $\bar{B}$, but only once the non-zero columns;*
- *all non-zero columns of $\bar{B}$ are present in $[\bar{B}_1^*, \ldots, \bar{B}_m^*]$;*

- *no two $\bar{B}_j^*, \bar{B}_k^*$, with $j,k = 1, \ldots, m$ and $j \neq k$, contain the same non-zero column vectors of $\bar{B}$,*

*then $\{\bar{B}_i^*\}_{i=1}^m$ is a solution to $\mathcal{P}_1$.*                                                          ⋄

**Proof:** By noticing that under the mentioned conditions the minimality is ensured, i.e., we have $\sum_{i=1}^m \|\bar{B}_i^*\|_0 = \|\bar{B}\|_0$, we only need to show the feasibility of the solution. Towards this goal, observe that there exists a $2mn \times 2mn$ permutation matrix $\mathbb{P}$ such that $[\bar{A}_1, \ldots, \bar{A}_m, \bar{B}_1^*, \ldots, \bar{B}_m^*]\mathbb{P} = [\bar{A}_1, \ldots, \bar{A}_m, \bar{B}, \mathbf{0}_{n \times n}, \ldots, \mathbf{0}_{n \times n}]$, where $\mathbf{0}_{n \times n}$ is the $n \times n$ zero matrix. Hence, by invoking Theorem 5 both conditions of Theorem 3 hold, and the result follows.                                                 ∎

The main results also provide new insights about the importance of the sequence of modes in the transition required to ensure structural controllability of the switching linear continuous-time systems. More specifically, we now do the following observations:

**Remark 2.** *Given a switching signal that ensures structural controllability of the switching linear continuous-time systems, the order of the modes among which the systems transitions does not impact the structural controllability of it. In fact, this follows from the conditions presented in Theorem 3, since both conditions are invariant to permutation.*                          ⋄

Subsequently, one may wonder which modes are the crucial to ensure both conditions in Theorem 3. This problem can be partially understood from the a solution obtained in Algorithm 1, which is captured in the following remark.

**Remark 3.** *Given the minimum weighted maximum matching $\mathcal{M}'$ obtained in Step 3 in Algorithm 1, if an edge $(c_\cdot^i, .) \in \mathcal{M}'$ then it follows that mode $i$ is being considered as part of the switching signal to ensure structural controllability of the switching linear continuous-time systems, since it contributes to ensure condition (ii) in Theorem 3. Nonetheless, this is not the same to say that there is no other minimum weighted maximum matching $\mathcal{M}''$ where $(c_\cdot^i, .) \notin \mathcal{M}''$. In other words, there might exist different switching signals ensuring structural controllability of the switching linear continuous-time systems, and these can be partially captured by the minimum weighted maximum matchings. Furthermore, it is possible to characterize which edges belong to any maximum matching [5], which implies that those modes need to be part of any sequence of modes among which the system transitions. Therefore, these modes should be considered as part of the design proposed in Theorem 3. In other words, the designer should consider to deploy actuation capabilities in the modes that are strictly required in a sequence to ensure structural* *controllability, i.e., at these modes, the input matrices should be non-zero.*                    ◇

Finally, we notice that determining the minimum sequence of modes a switching signal among which the system should transition to ensure structural controllability is NP-hard. Formally, consider the following problem.

$\mathcal{P}_2$ Given a structurally controllable $(\bar{A}_{\sigma(t)}, \bar{B}_{\sigma(t)})$, determine the minimum number of modes $m'$ that a switching signal $\sigma(t)$ should consider to attain structural controllability.                    ○

**Theorem 7.** *Problem $\mathcal{P}_2$ is NP-hard.*                    ◇

**Proof**: Consider the well known NP-hard problem, the set covering problem, that can be stated as follows: given a universe of elements $\mathcal{U} = \{1, \ldots, n\}$ and a collection of subsets $\{\mathcal{S}_j\}_{j=1}^m$ with $\mathcal{S}_j \subset \mathcal{U}$, determine the subcollection $\{\mathcal{S}_j\}_{j \in \mathcal{I}}$ that contains $\mathcal{U}$, where $\mathcal{I} \subset \{1, \ldots, m\}$ and there is no other $\mathcal{I}'$ such that $|\mathcal{I}'| < |\mathcal{I}|$ satisfying the same conditions.

In order to show that $\mathcal{P}_2$ is NP-hard, we need to polynomially reduce the set covering problem to $\mathcal{P}_2$ (see [6] for an introduction on the topic). As a consequence, a solution to $\mathcal{P}_2$ enables the reconstruction of a solution to the set covering problem, which implies that finding a solution to $\mathcal{P}_2$ is at least as difficult as finding a solution to the set covering problem. Towards this goal, we associate with each mode of the switching system a subset $\mathcal{S}_i$, and we assume that any mode suffices to ensure condition (ii) in Theorem 3; more specifically, we assume that $\bar{A}_i \in \{0, 1\}^{(n+1) \times (n+1)}$ and $\text{diag}(\bar{A}_i) = [1, \ldots, 1]_{(n+1) \times 1}$ for all $i = 1, \ldots, m$. Furthermore, let $\bar{A} = \bar{A}_1 \vee \ldots \vee \bar{A}_m$ be such that $\mathcal{D}(\bar{A}) = (\mathcal{V} = \{x_1, \ldots, x_{n+1}\})$ is a directed tree rooted in $x_1$; hence, by considering $[\bar{B}_1]_{1,1} = 1$ and $[\bar{B}_1]_{i,j} = 0$ for the remaining $i, j = 1, \ldots, n+1$, and $\bar{B}_k = 0$ for $k = 2, \ldots, m$, i.e., only the first mode is actuated and a dedicated input actuates $x_1$, it follows that condition (i) in Theorem 3 yields. Notice that there are $n$ directed edges in the directed tree, which we can enumerate as $\{e_i\}_{i=1}^n$. Therefore, each of the $m$ modes can exhibit in its digraph representation the edges $\{e_j\}_{j \in \mathcal{S}_k}$ for $k = 1, \ldots, m$, besides the self-loops in all state variables.

Under the present construction, it is not difficult to realize that any solution to $\mathcal{P}_2$ consists in finding the smallest subcollection of modes such that the system digraph $\mathcal{D}(\bar{A}, \bar{B})$, where $\bar{A} = \vee_{i \in \mathcal{I}^*} A_i$ and $\mathcal{I}^*$ denotes the indices of the modes that contain a directed spanning tree rooted in an input. Therefore, it follows that there exists a collection of edges producing such tree, which implies that there exist a collection of sets $\{\mathcal{S}_i\}_{i \in \mathcal{I}^*}$ that covers $\mathcal{U}$, and the result

follows.  ∎

## V. Illustrative Example

Consider a switched linear continuous-time system with three modes $\{\bar{A}_i\}_{i=1}^3$ where $\bar{A}_i \in \{0,1\}^{4\times 4}$, and $[A_1]_{1,2} = [A_2]_{3,2} = [A_3]_{4,4} = 1$ and zero otherwise. Then, we have

$$\bar{A} = \bar{A}_1 \vee \bar{A}_2 \vee \bar{A}_3 = \begin{bmatrix} 0 & 1 & 0 & 0 \\ 0 & 0 & 0 & 0 \\ 0 & 1 & 0 & 0 \\ 0 & 0 & 0 & 1 \end{bmatrix}.$$

Hereafter, we aim to determining the sparsest configuration of inputs that renders the system structurally controllable, i.e., a solution to $\mathcal{P}_1$. To this end, we consider Algorithm 1 (whose correctness and computational complexity are provided in Theorem 4).

First, the DAG representation of the state digraph associated with $\mathcal{D}(\bar{A})$ contains four SCCs depicted by dashed gray boxes in Figure 1-(d); in particular, $\mathcal{N}_1^T$ and $\mathcal{N}_2^T$ are two non-top linked SCC. Subsequently, we have $\bar{S}$ as follows:

$$\bar{S} = \begin{bmatrix} 0 & 0 \\ 1 & 0 \\ 0 & 0 \\ 0 & 1 \end{bmatrix}.$$

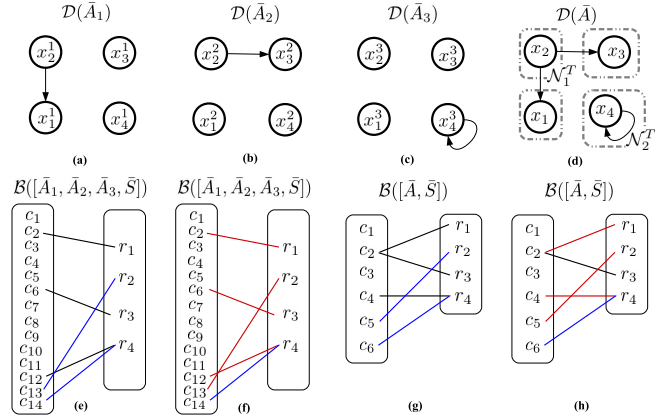

Fig. 1. In this figure, we depict the state digraph associated with $\bar{A}_i$, with $i = 1, 2, 3$, in (a)-(c), respectively. The state digraph of the union $\bar{A} = \bar{A}_1 \vee \bar{A}_2 \vee \bar{A}_3$ is depicted in (d) and contains three SCCs depicted by dashed gray boxes; in particular, $\mathcal{N}_1^T$ and $\mathcal{N}_2^T$ are two non-top linked SCC. In (e)-(f), we represent the bipartite graph obtained in Algorithm 1 and the MWMM. Alternatively, in (g)-(h) we illustrate the bipartite graph obtained in Algorithm 1 and the MWMM if a single mode is considered with $\bar{A}$ as its dynamics.

In addition, the weights associated with the edges in $\mathcal{B}([\bar{A}_1, \bar{A}_2, \bar{A}_3, \bar{S}])$ are as follows: the edges that contain the vertices $c_{13}$ and $c_{14}$ have unitary weight (depicted by the blue edges in Figure 1-(e)), and all other edges incur in zero weight (depicted by the black edges in Figure 1-(e)). At Step 3 in Algorithm 1, the MWMM is represented by the collection of red edges in Figure 1-(f), which we denote by $\mathcal{M}_1^*$. Lastly, taking $\mathcal{M}_1^*$, we obtain $\mathcal{J}' = \{2\}$, $\mathcal{J}'' = \emptyset$ (i.e., there is no need to increase the size of the maximum matching), and $\mathcal{J}''' = \{4\}$ that ensures the all nodes in $\mathcal{D}(\bar{A}, \mathbb{D}(\mathcal{J}' \cup \mathcal{J}'' \cup \mathcal{J}'''))$ are accessible.

Consequently, we obtain that $(\bar{B}_1, \bar{B}_2, \bar{B}_3) = (\mathbb{D}(\mathcal{J}' \cup \mathcal{J}'' \cup \mathcal{J}'''), 0, 0)$ is a dedicated solution to $\mathcal{P}_1$, by invoking Theorem 4. From Corollary 1, we obtain that $(\bar{B}_1', \bar{B}_2', \bar{B}_3') = (\mathbb{D}(\mathcal{J}' \cup \mathcal{J}'') \cup \mathbb{M}(\mathcal{J}'''), 0, 0)$ is a minimal solution to $\mathcal{P}_1$. If we want the actuation to be distributed among different modes of the switched linear continuous-time system, one just needs to recall Theorem 6; in particular, $(\bar{B}_1'', \bar{B}_2'', \bar{B}_3'') = (\mathbb{D}(\mathcal{J}'), \mathbb{D}(\mathcal{J}''), \mathbb{D}(\mathcal{J}''))$ is a solution to $\mathcal{P}_1$.

**Remark 4.** *The solution obtained for structural linear continuous-time switching systems cannot be retrieved from that of a structural linear time-invariant systems when a joint state digraph of $\bar{A}$ is considered. More specifically, consider a switching system with a single mode whose dynamics is $\bar{A}$, then by executing Algorithm 1, one needs to consider the bipartite graph in Figure 1-(g), where the edges in blue have unitary weight and the remaining have zero weight. Then, a possible MWMM determined in Step 3 in Algorithm 1 contains the edges depicted in red in Figure 1-(h), which implies that $\mathcal{J}' = \{2\}$, $\mathcal{J}'' = \{3\}$ and $\mathcal{J}''' = \{4\}$. In other words, an additional state variable is required to be actuated, i.e., $x_3$.* ◇

## VI. Conclusions and Further Research

In this brief paper, we provides a new necessary and sufficient condition that leverages both graph-theoretic and algebraic properties required to ensure structural controllability of switching linear continuous-time systems. With this condition we characterize the solutions to the structural minimum controllability problem for switched linear continuous-time systems. Further, we provided an efficient algorithm that determines a solution to the problem. Finally, we provides new insights on how the switching sequences affect the controllability of a structural switching linear continuous-time system. In particular, we show that determining the minimum collection of nodes in a sequence of modes ensuring structural controllability is NP-hard. Future research will focus on considering different actuation cost per state variables, actuators and possible switching sequences. Additionally, it would be interesting to address the sparsest feedback patterns that ensure the switched linear system to ensure stabilizability.